\documentclass[a4paper,11pt]{article}
\usepackage[utf8]{inputenc}

\usepackage{amsmath}
\usepackage{amsfonts}
\usepackage{amssymb}
\usepackage{enumerate}
\usepackage{graphics}
\usepackage{xcolor}
\RequirePackage{geometry}
\newtheorem{theorem}{Theorem}
\newtheorem{lemma}[theorem]{Lemma}

\newcommand{\prob}{\ensuremath{\mathbb{P}}}
\newcommand{\qed}{\hfill $\square$}

\newcommand{\RR}{\ensuremath{\mathbb{R}}}
\newcommand{\CC}{\ensuremath{\mathbb{C}}}
\newcommand{\NN}{\ensuremath{\mathbb{N}}}
\newcommand{\ZZ}{\ensuremath{\mathbb{Z}}}
\newcommand{\HP}{\ensuremath{\mathbb{H}}}
\newcommand{\TT}{\ensuremath{\mathbb{T}}}
\newcommand{\EE}{\mathbb{E}}
\newcommand{\lra}{\leftrightarrow}
\newcommand{\ind}{\mathbf{1}}

\title{The expected number of critical percolation clusters intersecting a line segment.}
\author{J. van den Berg\footnote{CWI and VU University Amsterdam, email: J.van.den.Berg@cwi.nl} \, and
R.P. Conijn\footnote{VU University Amsterdam (Current affiliation: Utrecht University, email: R.P.Conijn@uu.nl)}
}
\date{}

\begin{document}

\maketitle

\begin{abstract}
We study critical percolation on a regular planar lattice. Let $E_G(n)$ be the expected number of open clusters
intersecting or hitting the line segment $[0,n]$. (For the subscript $G$ we either take $\mathbb{H}$, when 
we restrict to
the upper halfplane, or $\mathbb{C}$, when we consider the full lattice). 

Cardy \cite{C01} (see also Yu, Saleur and Haas \cite{YSH08}) derived
heuristically that $E_{\mathbb{H}}(n) = An + \frac{\sqrt{3}}{4\pi}\log(n) + o(\log(n))$,
where $A$ is some constant.
Recently Kov\'{a}cs, Igl\'{o}i and Cardy derived in
\cite{KIC12} heuristically (as a special case of a more general formula) that a similar result holds for $E_{\mathbb{C}}(n)$
with the constant $\frac{\sqrt{3}}{4\pi}$ replaced by $\frac{5\sqrt{3}}{32\pi}$.

In this paper we give, for site percolation on the triangular lattice, a rigorous proof for the formula
of $E_{\mathbb{H}}(n)$ above, and a rigorous upper bound for the prefactor of the logarithm in the
formula of $E_{\mathbb{C}}(n)$.
\end{abstract}

\medskip\noindent
\begin{center}
{\small {\it 2010 Mathematics Subject Classification.} 60K35 (82B43). }
\end{center}

\section{Introduction}
\subsection{Background and statement of the main result}
Consider critical bond percolation on $\ZZ^2$. Kov\'{a}cs, Igl\'{o}i and Cardy \cite{KIC12} studied 
the expected number of clusters which intersect the boundary of a polygon. The leading order is
the size $n$ of the boundary. The prefactor of this term is lattice dependent. Their main interest is
in the first correction term (of order $\log n$). 
Their motivation came from relations with entanglement entropy in a diluted quantum Ising model.
Using indirect and non-rigorous methods from conformal field theory and the $q$-state Potts model (letting $q \to 1$),
they derived a (universal) formula for the prefactor of the logarithmic term.

A special case of their result is that of a line segment (treated in Section F of their paper). In their
setup the line segment was placed in the full plane
and they claim that the prefactor is equal to $\frac{5\sqrt{3}}{32\pi}$. Furthermore they refer
to an earlier obtained result by Cardy in \cite{C01} (see also Yu, Saleur and Haas \cite{YSH08})
where the line segment was placed on the boundary of the half-plane.
In the latter case the claim is that the prefactor equals $\frac{\sqrt{3}}{4\pi}$.
Also this latter result was obtained by non-rigorous arguments using $q$-state Potts models.

This motivated us to try to find rigorous and more direct proofs
of these results (starting with the case of line segments). Since the prefactors are believed to be
universal it is natural to consider the most well studied percolation model, site percolation on the
triangular lattice with $p = p_c = 1/2$.

Because conformal invariance plays a role, it is convenient to identify the plane with the set $\CC$ of complex numbers.
We embed the triangular lattice $\TT$ in the half-plane
$\HP = \{z: \Im z \ge 0\}$ or the full 
plane $\CC$ with vertex set $\{ m+n\mathbf{j}:\, m \in \ZZ, n\in \NN\cup\{0\}\}$
(resp. $\{ m+n\mathbf{j}:\, m,n \in \ZZ\}$), where $\mathbf{j} = e^{\frac{\pi}{3}\mathbf{i}}$.
We denote the probability measure by $\prob_{\HP}$
(resp. $\prob_{\CC}$) and the expectation by $\mathbb{E}_{\HP}$ (resp. $\mathbb{E}_{\CC}$).
For subsets $A,B \subset \CC$ we denote by $A\lra B$ the event that there are open vertices $x,y$
on the triangular lattice, with $x \in A, y \in B$, which are connected by a path of open vertices.
With some abuse of notation we denote,
for any $x \in \CC$, the set $\{x\}$ by $x$. A cluster is a maximal collection of connected vertices.
Consider the line segment $[1,n]$ on $\RR$, containing $n$ vertices.
We are interested in
\[
 E_{G}(n) := \EE_{G}[\,|\{C \in \mathcal{C}_{G}: C \cap [1,n] \neq \emptyset \}|\,],
\]
where $\mathcal{C}_{G}$ is the collection of clusters in the triangular lattice on the lattice $G = \HP, \CC$.

It is easy to derive the leading (of order $n$) term: see the Remark in Section \ref{subsec:IntroComputation}.
In the case of the half-plane we could obtain a rigorous proof for the earlier mentioned logarithmic correction term.
In the case of the full plane we only obtained a logarithmic upper bound for the correction term. (We do not see
a method how to
prove the precise prefactor $\frac{5\sqrt{3}}{32\pi}$ given in \cite{KIC12}; even finding a non-trivial lower bound
is, in our opinion, a challenging problem).

More precisely, our main contribution is a rigorous proof of the following:

\begin{theorem}\label{thm:mainResult}
 \begin{equation*}
  \textrm{(a)}\quad E_{\HP}(n) = n\cdot (\prob_{\HP}(1 \not\lra (-\infty,0])-\frac{1}{2}) + \frac{\sqrt{3}}{4\pi}\log(n) + o(\log(n))
 \end{equation*}
 and
 \begin{equation*}
\textrm{(b)}\qquad\,\,\,\,\,\,\,\,\,  \limsup_{n\to \infty} \frac{E_{\CC}(n) - n\cdot (\prob_{\CC}(1 \not\lra (-\infty,0])-\frac{1}{2})}{\log(n)} \le \frac{8}{5}\cdot \frac{\sqrt{3}}{4\pi}.
 \end{equation*}
\end{theorem}


\subsection{Some introductory computations}\label{subsec:IntroComputation}
We now describe the first steps of the strategy to derive the result above.
This will also give some insight, where the log comes from.
First rewrite the number of clusters as follows
\begin{eqnarray}
 |\{C \in \mathcal{C}_{G}: C \cap [1,n] \neq \emptyset \}| & = & \ind\{1 \textrm{ open}\} + \sum_{k=2}^{n} \ind\{k \not\lra [1,k-1], k \textrm{ open}\}\nonumber\\
 & = & 1 + \sum_{k=2}^{n} \,\ind\{k \not\lra [1,k-1]\} - \sum_{k=1}^{n} \ind\{k \textrm{ closed}\}\nonumber
\end{eqnarray}
So
\begin{eqnarray*}
 \lefteqn{E_{G}(n)}\\
 & = & 1 - \frac{1}{2}n + \sum_{k=2}^{n} \big(\prob_{G}(k \not\lra (-\infty,k-1]) + \prob_{G}(\,\{k \not\lra [1,k-1]\} \cap \{ k \lra (-\infty,0] \})\big)\\
      & = & 1 - \frac{1}{2}n + (n-1)\cdot (\prob_{G}(1 \not\lra (-\infty,0])) \\
      & & + \sum_{k=2}^{n} \prob_{G}(\,\{k \not\lra [1,k-1]\} \cap \{ k \lra (-\infty,0] \}).
\end{eqnarray*}

\smallskip\noindent
{\bf Remark:} It is known that there is no infinite cluster almost surely,
hence $ \prob_{G}(k \lra (-\infty,0]) \rightarrow 0$  as $k \rightarrow \infty$. Therefore the above computation implies that
the leading term of $E_{G}(n)$ is $n(\prob_{G}(1 \not\lra (-\infty,0]) - \frac{1}{2})$.

\smallskip
Let us introduce the following notation:
\[
 L_{G}(n) := \frac{1}{\log(n)}\sum_{k=2}^{n} \prob_{G}(\,\{k \not\lra [1,k-1]\} \cap \{ k \lra (-\infty,0] \}\,).
\]
That is,
\[
 L_{G}(n) =\frac{E_{G}(n) - 1 + \frac{1}{2}n - (n-1)\cdot (\prob_{G}(1 \not\lra (-\infty,0]))}{\log(n)}.
\]
Hence Theorem \ref{thm:mainResult} is equivalent to

\noindent (a) $\lim_{n\to \infty} L_{\HP}(n) = \frac{\sqrt{3}}{4\pi}$ and\newline
\noindent (b) $\limsup_{n \to \infty} L_{\CC}(n) \le \frac{8}{5}\cdot\frac{\sqrt{3}}{4\pi}$.

\smallskip
Take $\varepsilon > 0$. We will introduce $M = M(n,\varepsilon) \in \NN$ and a sequence $a(i) = a(i,n,\varepsilon)$ for $1 \le i \le M+1$,
such that \[ a(M+1) = n.\] With these values we split up the sum in $L_{G}(n)$ in the following terms.
For all $1 \le i \le M$,
\begin{equation}\label{eq:intro:def:fi}
 f_i := \sum_{k=a(i)+1}^{a(i+1)} \prob_{G}(\,\{k \not\lra [1,k-1]\} \cap \{ k \lra (-\infty,0] \}\,)
\end{equation}
and
\begin{equation}\label{eq:intro:def:f0}
 f_0 := \sum_{k=2}^{a(1)} \prob_{G}(\,\{k \not\lra [1,k-1]\} \cap \{ k \lra (-\infty,0] \}\,).
\end{equation}
Then 
\[
 L_{G}(n) = \frac{f_0}{\log(n)} + \frac{1}{\log(n)}\sum_{i=1}^{M} f_i.
\]
The idea is now, roughly speaking,
to choose $a(i,n,\varepsilon)$ so that the ratio of two consecutive ones equals $1+\varepsilon$
and choose $M$ such that $a(1,n,\varepsilon)$ goes to infinity as $n \to \infty$, but is of a smaller order than $\log(n)$.
Then obviously the term $f_0/\log(n)$ is negligible.
We will see that $M$ is more or less of the order $\log(n)/\varepsilon$.
The existence of the limit $\lim_{n\to \infty} L_{G}(n)$ would follow if we can show that,
for $\varepsilon$ close to zero, $f_i$ is approximately a constant times $\varepsilon$ as $n \to \infty$.

In the case that $G = \HP$, we will see in Section \ref{subsec:pf_HP} that this strategy indeed leads to the existence, and even the value,
of the limit of $L_{\HP}(n)$ as $n \to \infty$.
Unfortunately in the full-plane it only leads to the upper bound stated in Theorem \ref{thm:mainResult} (b),
as we will see in Section \ref{subsec:pf_CC}.

Now we make the above choices precise. We define
\begin{equation}\label{eq:def:M}
 M := \left\lfloor \frac{\log(n) - \frac{1}{2}\log(\log(n))}{\log(1+\varepsilon)}\right\rfloor
\end{equation}
and for $ i \in \{ -1,\cdots, M-1 \}$
\begin{equation}\label{eq:def:ai}
 a(M-i,n,\varepsilon) := \left\lfloor\frac{n}{(1+\varepsilon)^{i+1}}\right\rfloor 
\end{equation}
or alternatively, 
for $j \in \{ 1,\cdots, M+1 \}$
\[
 a(j,n,\varepsilon) := \left\lfloor\frac{n}{(1+\varepsilon)^{M-j+1}}\right\rfloor.
\]
Note that then $a(1,n,\varepsilon)$ is of order $\sqrt{\log(n)}$.
To examine $f_i$ it is useful to rewrite it in terms of an expectation as follows. Let
\begin{equation} \label{rob:eq-Tidef}
 T(i) := \sum_{k=a(i)+1}^{a(i+1)}\ind\{ k \not\lra [1,k-1] \textrm{ and } k \lra (-\infty,0] \}.
\end{equation}
Then $f_i = \EE_{G}[T(i)]$. Hence
\begin{equation}\label{eq:LGn-as-SumOf-Expec}
 L_{G}(n) = \frac{f_0}{\log(n)} + \frac{1}{\log(n)}\sum_{i=1}^{M} \EE_{G}[T(i)].
\end{equation}

\section{Ingredients from the literature}\label{sec:preliminaries}
In this section we state some results, which we will use in Section \ref{sec:PfMainResult} to prove Theorem \ref{thm:mainResult}.
First some additional notation.
We use the following notation for the probabilities of so-called arm-events. Let, for $m < n\in \NN$
\begin{equation}\label{eq:def:pi1}
 \pi_1(m,n) := \prob_{\HP}([-m,m]^2 \lra \HP\setminus [-n,n]^2)
\end{equation}
and let $\pi_3(m,n)$ be the probability of having two disjoint closed paths, and an open path, from
$[-m,m]^2$ to $\HP\setminus [-n,n]^2$.
The following lemma is well known (see for example Theorem 11, Proposition 14 and Theorem 24
 in \cite{N08}).
\begin{lemma}\label{lem:armBounds}
 There exist constants $C_1,C_2 > 0$ and $\alpha \le 1/2$ such that, for all $m < n$
 \[
  \pi_1(m,n) \le C_1\left(\frac{m}{n}\right)^{\alpha}, \qquad \pi_3(m,n) \le C_2\left(\frac{m}{n}\right)^2.
 \]
\end{lemma}
In fact, much more precise results for these probabilities are known, but will not be used in this paper.

In the rest of this section, for a simply connected domain $D \subsetneq \CC$ and $n \in \NN$ the notation $nD$ denotes the set
$\{ n\cdot u:\,u\in D\}$. For points $a_1,a_2$ on the boundary of $D$ we denote by $[a_1,a_2]$ the part
of the boundary of $D$ between $a_1$ and $a_2$ in the counter clockwise direction.
Furthermore we generalize the notation slightly, namely by $\prob_D$ (and $\EE_{D}$) we will denote
the probability measure for percolation restricted to the triangular lattice on $D$.
In this setting two intervals $[a_1,a_2]$ and $[a_3,a_4]$ on the boundary are said to be connected if there
are vertices $x,y$ on the lattice inside $D$, which are connected by an open path,
and are such that $x$ has an edge which crosses $[a_1,a_2]$ and $y$ has an edge which crosses $[a_3,a_4]$.

The first theorem is the famous Cardy's formula (proposed in \cite{C92}), which was proved by Smirnov in \cite{S01}.
\begin{theorem}[Cardy's formula, \cite{S01}]\label{thm:cardy} Let $D \subsetneq \CC$ be a simply connected domain
and $\phi: D \to \HP$ a conformal map. Let $a_1,a_2,a_3,a_4$ be ordered points on the boundary of $D$.
We have
\[
 \lim_{n\to \infty} \prob_{nD}([na_1,na_2] \lra [na_3,na_4]) = \frac{2\pi\sqrt{3}}{\Gamma\left(\frac{1}{3}\right)^3}\lambda^{1/3}\cdot \,_2F_1\left(\frac{1}{3},\frac{2}{3};\frac{4}{3};\lambda \right),
\]
where $\lambda$ is the cross-ratio
\begin{equation}\label{eq:lambdaCrossRatio}
 \lambda = \frac{(\phi(a_1)-\phi(a_2))(\phi(a_4)-\phi(a_3))}{(\phi(a_1)-\phi(a_3))(\phi(a_4)-\phi(a_2))}.
\end{equation}
\end{theorem}
This theorem concerns crossing probabilities of generalized rectangles in one 'direction'.
The following theorem gives a formula for probabilities of crossings in two directions.
It is called after Watts, who proposed the formula in \cite{W96}.
The first rigorous proof was by Dub\'{e}dat \cite{D06Watt}. An alternative proof was obtained by Schramm (see
\cite{SW11}).
\begin{theorem}[Watts' formula, \cite{D06Watt,SW11}]\label{thm:watts} Let $D \subsetneq \CC$ be a simply
connected domain and $\phi: D \to \HP$ a conformal map. Let $a_1,a_2,a_3,a_4$ be ordered points on
the boundary of $D$. We have
\begin{eqnarray*}
 \lefteqn{\lim_{n\to \infty} \prob_{nD}([na_1,na_2] \lra [na_3,na_4]\textrm{ and } [na_4,na_1] \lra [na_2,na_3])} \\
 & = & \frac{2\pi\sqrt{3}}{\Gamma\left(\frac{1}{3}\right)^3}\lambda^{1/3}\cdot \,_2F_1\left(\frac{1}{3},\frac{2}{3};\frac{4}{3};\lambda \right) - \frac{\sqrt{3}}{2\pi}\lambda\cdot \,_3F_2\left(1,1,\frac{4}{3};\frac{5}{3},2;\lambda\right),
\end{eqnarray*}
 where $\lambda$ is the cross-ratio \eqref{eq:lambdaCrossRatio}. 
\end{theorem}
The last theorem we state here concerns the expected number of crossing clusters of a rectangle.
It was predicted by Cardy \cite{C01} and by Simmons, Kleban and Ziff \cite{SKZ07b}.
A proof was given by Hongler and Smirnov in \cite{HS11}.
Here $N(nD,a_1,a_2,a_3,a_4)$ denotes the number of clusters in $nD$ which connect $[na_1,na_2]$ with $[na_3,na_4]$.
\begin{theorem}[\cite{HS11}]\label{thm:ExpCrossClusters} Let $D \subsetneq \CC$ be a simply connected domain
and $\phi: D \to \HP$ a conformal map. Let $a_1,a_2,a_3,a_4$ be ordered points on the boundary of $D$.
We have
\begin{eqnarray*}
 \lefteqn{\lim_{n\to \infty} \EE_{nD}[N(nD,a_1,a_2,a_3,a_4)]} \\
 & = & \frac{2\pi\sqrt{3}}{\Gamma\left(\frac{1}{3}\right)^3}\lambda^{1/3}\cdot \,_2F_1\left(\frac{1}{3},\frac{2}{3};\frac{4}{3};\lambda \right) - \frac{\sqrt{3}}{4\pi}\lambda\cdot \,_3F_2\left(1,1,\frac{4}{3};\frac{5}{3},2;\lambda\right) + \frac{\sqrt{3}}{4\pi}\log\left(\frac{1}{1-\lambda}\right),
\end{eqnarray*}
 where $\lambda$ is the cross-ratio \eqref{eq:lambdaCrossRatio}. 
\end{theorem}

\section{Proof of Theorem \ref{thm:mainResult}}\label{sec:PfMainResult}
Recall from the introduction that Theorem \ref{thm:mainResult} is equivalent to

\noindent (a) $\lim_{n\to \infty} L_{\HP}(n) = \frac{\sqrt{3}}{4\pi}$ and\newline
\noindent (b) $\limsup_{n \to \infty} L_{\CC}(n) \le \frac{8}{5}\cdot\frac{\sqrt{3}}{4\pi}$.

Recall the definition \eqref{rob:eq-Tidef} of $T(i)$.
We begin this section with a lemma which says that, to prove
the convergence of $L_{G}(n)$ as $n\to \infty$, it is sufficient to prove
the convergence of $\varepsilon^{-1}\EE_{G}[T(i)]$.

\begin{lemma}\label{lem:sufficient} The following inequalities hold.
 \begin{equation}\label{eq:lem:Part1}
  \limsup_{n\to \infty} L_{G}(n) \le \limsup_{\varepsilon \to 0}\, \limsup_{n\to \infty}\, \max_{1 \le i \le M} \frac{\EE_{G}[T(i)]}{\varepsilon}
 \end{equation}
 and
 \begin{equation}\label{eq:lem:Part2}
  \liminf_{n\to \infty} L_{G}(n) \ge \liminf_{\varepsilon \to 0}\, \liminf_{n\to \infty}\, \min_{1 \le i \le M} \frac{\EE_{G}[T(i)]}{\varepsilon}.
 \end{equation}
\end{lemma}

\noindent\textit{Proof: }Recall \eqref{eq:LGn-as-SumOf-Expec} and the definitions of $M,a(i),f_i$ in
\eqref{eq:intro:def:fi} - \eqref{eq:def:ai}.
To prove \eqref{eq:lem:Part1}, first note that
$0 \le f_0 \le a(1,n,\varepsilon)$ and $M$ was chosen such that $a(1,n,\varepsilon)\approx \sqrt{\log(n)}$,
hence
\[
 \lim_{n\to \infty}\, \frac{f_0}{\log(n)} = 0.
\]
Thus it is enough to prove that
\begin{equation}\label{eq:pf:lem:enough}
 \limsup_{\varepsilon\to 0}\, \limsup_{n\to \infty}\, \left( \sum_{i=1}^{M} \frac{\EE_{G}[T(i)]}{\log(n)} \right) \le 
  \limsup_{\varepsilon \to 0}\, \limsup_{n\to \infty}\, \max_{1 \le i \le M} \frac{\EE_{G}[T(i)]}{\varepsilon}.
\end{equation}
Hereto, note that it is also easy to see from the definition of $M$ that, for fixed $\varepsilon > 0$
\[
 \lim_{n\to \infty} \frac{M}{\log(n)} = \frac{1}{\log(1+\varepsilon)}.
\]
For all $\varepsilon > 0$ we have
\begin{eqnarray}
\limsup_{n \to \infty} \sum_{i=1}^{M} \frac{\EE_{G}[T(i)]}{\log(n)} & \le & \limsup_{n\to \infty}\left( \frac{M}{\log(n)}\max_{i\le M} \EE_{G}[T(i)]\right) \nonumber\\
       & \le & \frac{1}{\log(1+\varepsilon)}\cdot \varepsilon\cdot \limsup_{n\to\infty}\left( \max_{i\le M} \frac{\EE_{G}[T(i)]}{\varepsilon}\right).\label{eq:pf:lem:limsupn}
\end{eqnarray}
Next note that
\[
 \limsup_{\varepsilon \to 0} \left( \frac{\varepsilon}{\log(1+\varepsilon)}\cdot \limsup_{n\to\infty}\left( \max_{i\le M} \frac{\EE_{G}[T(i)]}{\varepsilon}\right)\right)
  =  \limsup_{\varepsilon \to 0} \limsup_{n\to\infty}\left( \max_{i\le M} \frac{\EE_{G}[T(i)]}{\varepsilon}\right).
\]
This together with \eqref{eq:pf:lem:limsupn} implies \eqref{eq:pf:lem:enough} and completes the proof of
\eqref{eq:lem:Part1}.


The inequality in \eqref{eq:lem:Part2} follows in a similar way and we omit it.\qed

\subsection{Proof of Theorem \ref{thm:mainResult} (a)}\label{subsec:pf_HP}
\begin{figure}\label{fig:Tge2}
 \centering
 \scalebox{0.6}{\includegraphics{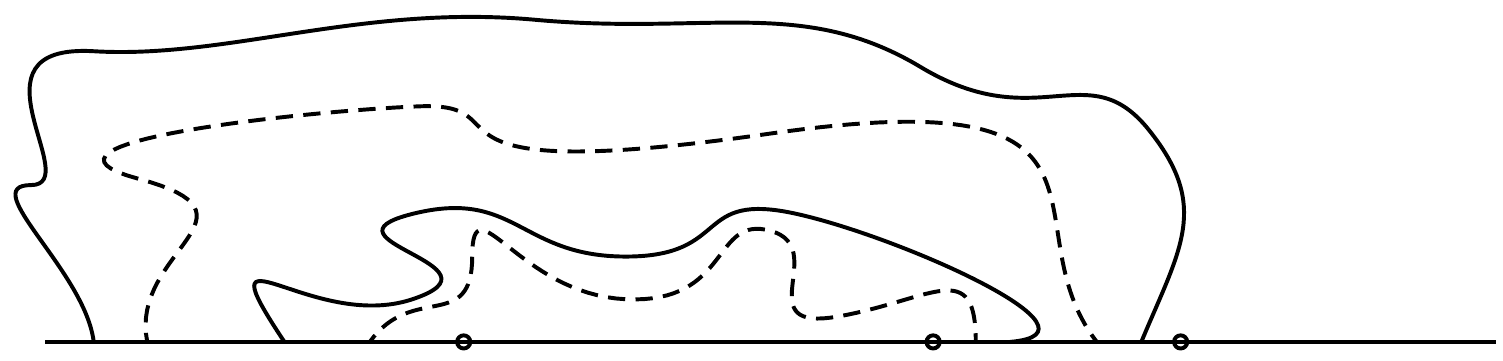}}
  \caption{The event $\{T(i) \ge 2\}$ occurs.}
\end{figure}
First note that it is easy to see that
$\{T(i) \ge 1\}$ if and only if there is an open and a closed path
from $(-\infty,1]$ to $[a(i),a(i+1)]$
and the closed path is below the open path.
Furthermore the event $\{T(i) \ge m\}$ is equal to the event that there are $2m$ alternating paths
between the aforementioned intervals, starting, from below, with a closed path. See Figure 1.
By conditioning on the lowest closed path and the lowest open path above that path, it is easy to see that, for each m,
$\prob_{\HP}(T(i) \geq m) \leq \prob_{\HP}(T(i) \geq 1) \times \prob_{\HP}(T(i) \geq m-1)$, and hence
\begin{equation}\label{eq:pf:Thm:a:BKusage}
 \prob_{\HP}(T(i) \ge m) \le (\prob_{\HP}(T(i) \ge 1))^m.
\end{equation}
Therefore
\begin{eqnarray}
 \EE_{\HP}[T(i)] & = & \sum_{m=1}^{\infty} \prob_{\HP}(T(i) \ge m) \nonumber\\
 & \le & \prob_{\HP}(T(i) \ge 1) + \sum_{m=2}^{\infty} (\prob_{\HP}(T(i) \ge 1))^m \nonumber\\
 & = & \prob_{\HP}(T(i) \ge 1) + \frac{(\prob_{\HP}(T(i) \ge 1))^2}{1 - \prob_{\HP}(T(i) \ge 1)}.\label{eq:pf:Thm:a:expect}
\end{eqnarray}

It is well-known from standard RSW arguments that $\prob_{\HP}(T(i) \ge 1)$ goes, uniformly in $i$ and $n$, to $0$
as $\epsilon \rightarrow 0$, since the ratio between two consecutive $a(i)$'s goes
to 1 as $\varepsilon \to 0$. Hence the `error term' (i.e. the second term in the r.h.s. of the equation array above)
is negligible w.r.t. the main term (i.e. the first term in the r.h.s.). By this, Lemma \ref{lem:sufficient},
the fact that $a(1) \rightarrow \infty$ as $n \rightarrow\infty$, and the ratio between consecutive $a(i)$'s,
it is sufficient to prove that

\begin{equation} \label{rob:eq-wk}
\lim_{k \rightarrow \infty} \prob_{\HP}(W_k) = \frac{\sqrt{3}}{4\pi} \varepsilon + o(\varepsilon),
\end{equation}
where $W_k$ denotes the event that there is an open and a closed path from 
$(-\infty,1]$ to $[k, k(1 + \varepsilon)]$ and the closed path is below the open path.

Let $W'_k$ be the event that there is an open and a closed path from
$(-\infty,1]$ to $[k, k(1 + \varepsilon)]$. (So, informally speaking, $W'_k$ is the same as $W_k$ without the
condition on which path is above or below).
Using that (by duality), there is either an open path from $[1,k]$ to $[k (1 + \varepsilon), \infty)$ or a closed
path from $(-\infty,1]$ to $[k, k(1 + \varepsilon)]$, we have

\begin{eqnarray} \label{rob:eq-duality}
 \lefteqn{\prob_{\HP}((-\infty,0] \lra [k, k(1 + \varepsilon)] \mbox{ and } [1,k] \lra [k(1 + \varepsilon), \infty))}
\;\;\;\;\;\;\quad\qquad\qquad \nonumber\\
 & = &  \prob_{\HP}((-\infty,1] \lra [k, k(1 + \varepsilon)]) - \prob_{\HP}(W'_k).
\end{eqnarray}
The limits as $k \rightarrow \infty$ of the first probability in the r.h.s. and the probability in the l.h.s.
are obtained by Theorem \ref{thm:cardy} and Theorem \ref{thm:watts} respectively, and we get

\begin{eqnarray} \label{rob:eq-w'k}
\lim_{k \rightarrow \infty} \prob_{\HP}(W'_k)
& = & \frac{\sqrt{3}}{2\pi}\cdot\frac{\varepsilon}{1+\varepsilon}\cdot
\,_3F_2\left(1,1,\frac{4}{3};\frac{5}{3},2;\frac{\varepsilon}{1+\varepsilon}\right) \nonumber\\
 & = & 2 \frac{\sqrt{3}}{4\pi}\cdot \varepsilon + o(\varepsilon).
\end{eqnarray}

Finally, let $\tilde W_k$ denote the event obtained from $W_k$ by replacing `open' by `closed' and vice versa.
Since $W_k$ and $\tilde W_k$ have the same probability and $W'_k = \tilde{W}_k \cup W_k$,
we have 
\begin{equation}\label{eq:EqualityWk}
 \prob_{\HP}(W'_k) = 2 \prob_{\HP}(W_k) - \prob_{\HP}(W_k \cap \tilde W_k).
\end{equation}
Since $W_k \cap \tilde W_k$ is contained in the disjoint occurrence of $W'_k$ and the event that there is
an open or closed path from $(-\infty,1]$ to $[k, k(1 + \varepsilon)]$, its probability is negligible
(as $k \to \infty$ and $\varepsilon \to 0$) w.r.t. that of $W'_k$, and we get from \eqref{rob:eq-w'k} and \eqref{eq:EqualityWk} that
$$\lim_{k \rightarrow \infty} \prob_{\HP}(W_k) = \frac{\sqrt{3}}{4\pi}\cdot \varepsilon +o(\varepsilon).$$
As we saw (see the argument above \eqref{rob:eq-wk}) this proves
Theorem \ref{thm:mainResult} (a). \qed

\subsection{Proof of Theorem \ref{thm:mainResult} (b)}\label{subsec:pf_CC}
We will bound the relevant probabilities (concerning the full plane) by the probabilities of certain
connection events in the half-plane. We do this by cutting along the real line from $-\infty$ up to $a(i+1)$.
Let us make the cutting precise.
Let
\[
 L(i) := (-\infty,a(i+1)],
\]
we define the new lattice to be the triangular lattice on $\CC \setminus L(i)$. This is the full triangular
lattice, without the vertices (and their edges) on $L(i)$.
Let us denote the corresponding probability measure, concerning percolation on this sublattice,
by $\tilde{\prob_{i}}$
(and expectation by $\tilde{\mathbb{E}_{i}}$).
Let the boundary $\partial_{\TT}[a,b]$ of an interval $[a,b] \subset L(i)$ be
the vertices $v$ of $\TT$ which are not in the interval $[a,b]$ but have a neighbouring vertex
which is on the interval $[a,b]$.
Let $\tilde{T}(i)$ be the number of clusters which connect $\partial_{\TT} [a(i)+1,a(i+1)]$
with $\partial_{\TT} (-\infty,0]$ but are not connected with $\partial_{\TT} [1,a(i)]$.

\begin{figure}\label{fig:probBi}
 \centering
 \scalebox{1.2}{\includegraphics{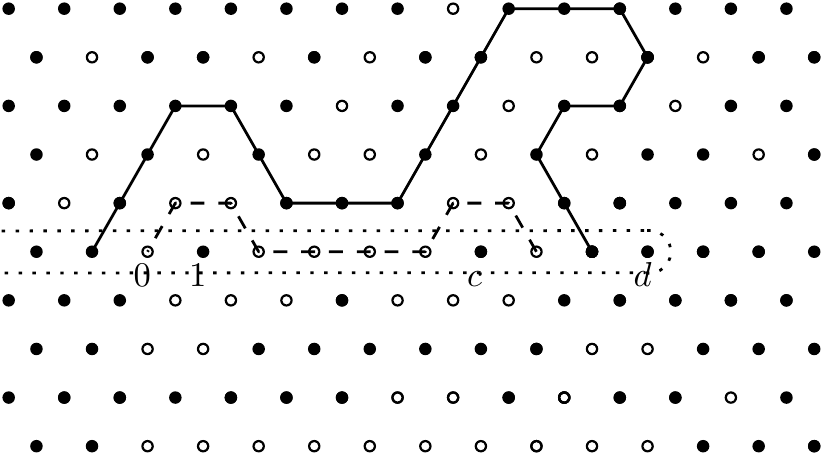}}
  \caption{The event $B(i)$ occurs. Filled circles are open and empty circles are closed.
  $c = a(i), d=a(i+1)$.}
\end{figure}
With this definition of $\tilde{T}(i)$ `almost all' the open connections counted in $T(i)$
are counted in $\tilde{T}(i)$ as well; however, there are exceptions.
In these exceptional cases there is an open connection from $(-\infty, 0]$ to $[a(i)+1, a(i+1)]$
which is not connected to $[1,a(i)]$ on $\TT$, but \'{i}s connected
to $\partial_{\TT} [1,a(i)]$ on $\CC \setminus L(i)\cap \TT$. See Figure 2.
More precisely, we define
\[
 B(i) := \bigcup_{k \in [1,a(i)] \cap \TT} \left(B_u(i,k) \cup B_l(i,k)\right),
\]
where $B_u(i,k)$ is the event that, on $\HP \cap \TT$,
there are closed paths from $k$ to $(-\infty, 1]$ and from $k$ to $[a(i),\infty)$ and open paths from one of the vertices
$k+\mathbf{j}$ and $k-1+\mathbf{j}$ to $(-\infty, 0]$ and to $[a(i)+1, \infty)$. (The open paths are not necessarily disjoint).
The event $B_l(i,k)$ is defined similarly on the lower half-plane.

We have
\begin{equation}\label{eq:pf:CC:ETTtilde}
 \EE_{\CC}[T(i)] \le \tilde{\mathbb{E}_{i}}[\tilde{T}(i)] + 2\prob_{\CC}(B(i)).
\end{equation}
To bound $\prob_{\CC}(B(i))$ we use the first inequality of Lemma \ref{lem:armBounds} for those $k$ in the
definition of $B(i)$ that are `close to' $1$ or
$a(i)$, and the other inequality in that lemma for the other $k$'s. More precisely, we fix a constant $\beta \in (0,1)$,
and let $r(a(i)) := \lceil a(i)^\beta\rceil$.
Then,
\begin{eqnarray}
 \prob_{\CC}(B(i)) & \le & 4\pi_1(r(a(i)),a(i)) +4\sum_{k=r(a(i)) + 1}^{\lceil\frac{1}{2}a(i)\rceil} \pi_3(1,k)\nonumber\\
 & \le & 4C_1\left(\frac{r(a(i))}{a(i)}\right)^{\alpha} +4\sum_{k=r(a(i)) + 1}^{\infty}  C_2\left(\frac{1}{k}\right)^2,\label{eq:pf:partb:boundB}
\end{eqnarray}
where the factor $4$ comes from symmetry considerations.
Hence, there exist constants $C_3, C_4 > 0$ such that
\begin{equation}\label{eq:pf:partb:boundB2}
 \prob_{\CC}(B(i)) \le C_3\left(a(i)^{\beta-1}\right)^{\alpha} + \frac{C_4}{a(i)^\beta}.
\end{equation}
Note that, since $a(1)$ (the smallest of the $a(i)$'s) tends to $\infty$ as $n \rightarrow \infty$, and
$C_3 (x^{\beta - 1})^{\alpha} + \frac{C_4}{x^\beta}$ tends to $0$ as $x \rightarrow \infty$, the
contribution of $\prob_{\CC}(B(i))$ to the r.h.s. of \eqref{eq:lem:Part1} is $0$.

\begin{figure}\label{fig:RelTandS}
 \centering
 \scalebox{1.2}{\includegraphics{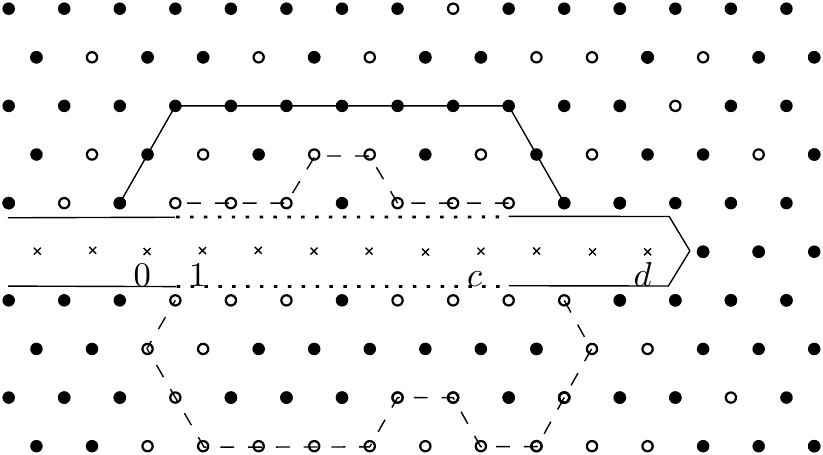}}
  \caption{Illustration of the event $\tilde{T}(i)=1$; here $c = a(i), d=a(i+1)$. The half-line $L(i)$ cut out of
  the lattice is indicated by x's. The boundaries
$\partial_{\TT} (-\infty,0]$ and $\partial_{\TT} [a(i)+1,a(i+1)]$ are indicated by solid lines and the boundary $\partial_{\TT} [1,a(i)]$ by dotted lines.
Filled circles are open and empty circles are closed. The solid path denotes an open connection and the dashed paths denote closed connections.}
\end{figure}
Next we consider the term $\tilde{\mathbb{E}_{i}}[\tilde{T}(i)]$.
Let $S(i)$ denote the number of closed clusters (on the earlier mentioned sublattice on $\mathbb{C}\setminus L(i)$)
connecting $\partial_{\TT} [a(i)+1,a(i+1)]$ with $\partial_{\TT} (-\infty,0]$.
Note that $\partial_{\TT}[1,a(i)]$ consists of two separate pieces (one in the upper and one in the lower half-plane) and observe that if there is only one open cluster as in the definition of $\tilde{T}(i)$, there are two closed clusters as in the definition of $S(i)$ (one preventing the open cluster to touch the mentioned `upper piece' of $\partial_{\TT}[1,a(i)]$, and one preventing it to touch the `lower piece'; see e.g. Figure 3).  So in this case we have $S(i) =2$.
More generally, a similar observation gives
\[
 \tilde{T}(i) = S(i) - \ind\{S(i) \ge 1\}.
\]
Thus it follows immediately, that
\begin{equation}\label{eq:pf:CC:TS}
 \tilde{\mathbb{E}_{i}}[\tilde{T}(i)] = \tilde{\mathbb{E}_{i}}[S(i)] - \tilde{\prob_{i}}(S(i) \ge 1).
\end{equation}
To complete the proof we will use Theorem \ref{thm:ExpCrossClusters}. Therefore we consider
the domain $\CC \setminus L(i)$ and scale it by $a(i)$. (As noted before, $a(1)$ goes to $\infty$
as $n \rightarrow \infty$). This gives the conformal rectangle
$\CC\setminus (-\infty,1+\varepsilon)$ with `corners' $a_1 = 1^+$, $a_2 = 1^-$, $a_3 = 0^-$, $a_4 = 0^+$ (where, for $x < 1 + \varepsilon$,
$x^+$ and $x^-$ denote the `copy' of $x$ in the upper and the lower half-plane respectively).
To apply Theorem \ref{thm:ExpCrossClusters} we need the cross-ratio, which can be computed as follows:
Consider
the conformal map
\[
 \varphi(z) := \mathbf{i}\sqrt{z-1-\varepsilon}
\]
which maps $\CC\setminus (-\infty,1+\varepsilon)$ onto the upper half-plane.
The cross-ratio is 
\[
 \lambda(\varepsilon) = \frac{(\varphi(1^+) - \varphi(1^-))(\varphi(0^+) - \varphi(0^-))}{(\varphi(1^+) - \varphi(0^-))(\varphi(0^+) - \varphi(1^-))}.
\]
It is easy to see that
\[
 \varphi(0^-) = \sqrt{1+\varepsilon},\quad \varphi(1^-) = \sqrt{\varepsilon},\quad \varphi(1^+) = -\sqrt{\varepsilon},\quad \varphi(0^+) = -\sqrt{1+\varepsilon}.
\]
Hence
\begin{eqnarray}
\lambda(\varepsilon)^2 & = & \frac{16\varepsilon(1+\varepsilon)}{(\sqrt{1+\varepsilon}+\sqrt{\varepsilon})^4} \nonumber\\
                       & = & 16\varepsilon + o(\varepsilon).\label{eq:pf:CC:lambdasquared}
\end{eqnarray}

Applying Theorem \ref{thm:ExpCrossClusters} we conclude that, as $n \rightarrow\infty$, $\tilde{\mathbb{E}_{i}}[S(i)]$
converges (uniformly in the $i's$ with $1 \leq i \leq M(n)$), to
\begin{eqnarray*}
\lefteqn{ \frac{2\pi\sqrt{3}}{\Gamma(\frac{1}{3})^3} \lambda(\varepsilon)^{1/3}\cdot \,_2F_1\left(\frac{1}{3},\frac{2}{3};\frac{4}{3};\lambda(\varepsilon)\right)} \\
 & - & \frac{\sqrt{3}}{4\pi}\lambda(\varepsilon)\cdot \,_3F_2\left(1,1,\frac{4}{3};\frac{5}{3},2;\lambda(\varepsilon)\right)+\frac{\sqrt{3}}{4\pi}\log\left(\frac{1}{1-\lambda(\varepsilon)}\right).
\end{eqnarray*}
The first term is exactly the limit $\tilde{\prob_{i}}(S(i) \ge 1)$ as $n \to \infty$ (Cardy's formula).
Hence by noting that
\[
  - \frac{\sqrt{3}}{4\pi}\lambda\cdot \,_3F_2\left(1,1,\frac{4}{3};\frac{5}{3},2;\lambda\right)+\frac{\sqrt{3}}{4\pi}\log\left(\frac{1}{1-\lambda}\right) = \frac{\sqrt{3}}{4\pi}\cdot \frac{1}{10}\lambda^2 + o(\lambda^2),
\]
and \eqref{eq:pf:CC:TS} and \eqref{eq:pf:CC:lambdasquared} we get that
\begin{equation}\label{eq:pf:CC:atLeasttwoCrossings}
  \lim_{n\to \infty} \tilde{\mathbb{E}_{i}}[\tilde{T}(i)] = 
\frac{\sqrt{3}}{4\pi}\cdot \frac{16}{10}\varepsilon + o(\varepsilon) =
\frac{8}{5}\cdot\frac{\sqrt{3}}{4\pi}\cdot \varepsilon + o(\varepsilon),
\end{equation}
uniformly in the $i$'s with $1 \leq i \leq M(n)$.

This, combined with \eqref{eq:pf:CC:ETTtilde} and the negligibility of $\prob_{\CC}(B(i))$
(see the line below \eqref{eq:pf:partb:boundB2}), gives

\[
 \limsup_{n \to \infty}\,\max_{1 \le i \le M} \EE_{\CC}[T(i)] \le \frac{8}{5}\cdot\frac{\sqrt{3}}{4\pi}\cdot \varepsilon + o(\varepsilon).
\]
By Lemma \ref{lem:sufficient} this implies Theorem \ref{thm:mainResult} (b).\qed

\medskip\noindent
{\bf \large Acknowledgment.} The first author thanks Christian Maes
for drawing his attention to the paper by Kov\'{a}cs, Igl\'{o}i and Cardy, in the spring of 2013.

\bibliographystyle{amsalpha}
\bibliography{mybiblio}

\providecommand{\bysame}{\leavevmode\hbox to3em{\hrulefill}\thinspace}
\providecommand{\MR}{\relax\ifhmode\unskip\space\fi MR }
\providecommand{\MRhref}[2]{%
  \href{http://www.ams.org/mathscinet-getitem?mr=#1}{#2}
}
\providecommand{\href}[2]{#2}
\begin{thebibliography}{YSH08}

\bibitem[Car92]{C92}
J.L. Cardy, \emph{Critical percolation in finite geometries}, Journal of
  Physics A: Mathematical and General \textbf{25} (1992), no.~4, L201--L206.

\bibitem[Car01]{C01}
\bysame, \emph{Lectures on conformal invariance and percolation}, Arxiv:0103018
  (2001).

\bibitem[Dub06]{D06Watt}
Julien Dub{\'e}dat, \emph{Excursion decompositions for {SLE} and {W}atts'
  crossing formula}, Probab. Theory Related Fields \textbf{134} (2006), no.~3,
  453--488. \MR{2226888 (2007d:60019)}

\bibitem[HS11]{HS11}
Cl{\'e}ment Hongler and Stanislav Smirnov, \emph{Critical percolation: the
  expected number of clusters in a rectangle}, Probability theory and related
  fields \textbf{151} (2011), no.~3-4, 735--756.

\bibitem[KIC12]{KIC12}
Istv\'{a}n~A. Kov\'{a}cs, Ferenc Igl\'{o}i, and John Cardy, \emph{Corner
  contribution to percolation cluster numbers}, Phys. Rev. B \textbf{86}
  (2012), 1--6.

\bibitem[Nol08]{N08}
P.~Nolin, \emph{Near-critical percolation in two dimensions}, Electron. J.
  Probab. \textbf{13} (2008), no. 55, 1562--1623. \MR{2438816 (2009k:60215)}

\bibitem[SKZ07]{SKZ07b}
Jacob~JH Simmons, Peter Kleban, and Robert~M Ziff, \emph{Percolation crossing
  formulae and conformal field theory}, Journal of Physics A: Mathematical and
  Theoretical \textbf{40} (2007), no.~31, F771.

\bibitem[Smi01]{S01}
Stanislav Smirnov, \emph{Critical percolation in the plane: conformal
  invariance, {C}ardy's formula, scaling limits}, C. R. Acad. Sci. Paris S\'er.
  I Math. \textbf{333} (2001), no.~3, 239--244. \MR{1851632 (2002f:60193)}

\bibitem[SW11]{SW11}
Scott Sheffield and David~B Wilson, \emph{Schramm's proof of {W}atts' formula},
  The Annals of Probability \textbf{39} (2011), no.~5, 1844--1863.

\bibitem[Wat96]{W96}
G.M.T. Watts, \emph{A crossing probability for critical percolation in two
  dimensions}, Journal of Physics A: Mathematical and General \textbf{29}
  (1996), no.~14, L363.

\bibitem[YSH08]{YSH08}
Rong Yu, Hubert Saleur, and Stephan Haas, \emph{Entanglement entropy in the
  two-dimensional random transverse field {I}sing model}, Physical Review B
  \textbf{77} (2008), no.~14, 140402.

\end{thebibliography}

\end{document}